\documentclass[10pt,letterpaper]{article}
\usepackage[top=0.85in,left=2.75in,footskip=0.75in]{geometry}

\usepackage{amsmath,amssymb}

\usepackage{changepage}

\usepackage[utf8x]{inputenc}

\usepackage{textcomp,marvosym}

\usepackage{cite}

\usepackage{nameref,hyperref}

\usepackage[right]{lineno}

\usepackage{microtype}
\DisableLigatures[f]{encoding = *, family = * }

\usepackage[table]{xcolor}

\usepackage{array}

\newcolumntype{+}{!{\vrule width 2pt}}

\newlength\savedwidth

\newcommand\thickhline{\noalign{\global\savedwidth\arrayrulewidth\global\arrayrulewidth 2pt}%
\hline
\noalign{\global\arrayrulewidth\savedwidth}}

\raggedright
\usepackage[aboveskip=1pt,labelfont=bf,labelsep=period,justification=raggedright,singlelinecheck=off]{caption}

\makeatletter
\renewcommand{\@biblabel}[1]{\quad#1.}
\makeatother

\usepackage{lastpage,fancyhdr,graphicx}
\usepackage{epstopdf}
\fancyheadoffset[L]{2.25in}
\fancyfootoffset[L]{2.25in}
\nolinenumbers

\begin{document}
\vspace*{0.2in}

\begin{flushleft}
{\Large
\textbf\newline{Comparing demographics of signatories to public letters on diversity in the mathematical sciences} 
}
\newline
\\
Chad M. Topaz \textsuperscript{1,2*},
James Cart\textsuperscript{3},
Carrie Diaz Eaton\textsuperscript{4},
Anelise Hanson Shrout\textsuperscript{4},
Jude A. Higdon\textsuperscript{2},
Kenan {\.{I}}nce\textsuperscript{5},
Brian Katz\textsuperscript{6},
Drew Lewis\textsuperscript{7},
Jessica Libertini\textsuperscript{8},
Christian Michael Smith\textsuperscript{9}
\\
\bigskip
\textbf{1} Department of Mathematics and Statistics, Williams College, Williamstown, MA, USA
\\
\textbf{2} Institute for the Quantitative Study of Inclusion, Diversity, and Equity, Williamstown, MA, USA
\\
\textbf{3} Office of Institutional Research, Williams College, Williamstown, MA, USA
\\
\textbf{4} Digital and Computational Studies Program, Bates College, Lewiston, ME, USA
\\
\textbf{5} Department of Mathematics, Westminster College, Salt Lake City, UT, USA
\\
\textbf{6} Department of Mathematics and Statistics, Smith College, Northampton, MA, USA
\\
\textbf{7} Department of Mathematics and Statistics, University of South Alabama, Mobile, AL, USA
\\
\textbf{8} Department of Applied Mathematics, Virginia Military Institute, Lexington, VA, USA
\\
\textbf{9} Department of Sociology, University of Wisconsin-Madison, Madison, WI, USA
\\
\bigskip

* qside@qsideinstitute.org

\end{flushleft}
\section*{Abstract}
In its December 2019 edition, the \textit{Notices of the American Mathematical Society} published an essay critical of the use of diversity statements in academic hiring. The publication of this essay prompted many responses, including three public letters circulated within the mathematical sciences community. Each letter was signed by hundreds of people and was published online, also by the American Mathematical Society. We report on a study of the signatories' demographics, which we infer using a crowdsourcing approach. Letter A highlights diversity and social justice. The pool of signatories contains relatively more individuals inferred to be women and/or members of underrepresented ethnic groups. Moreover, this pool is diverse with respect to the levels of professional security and types of academic institutions represented. Letter B does not comment on diversity, but rather, asks for discussion and debate. This letter was signed by a strong majority of individuals inferred to be white men in professionally secure positions at highly research intensive universities. Letter C speaks out specifically against diversity statements, calling them ``a mistake,'' and claiming that their usage during early stages of faculty hiring ``diminishes mathematical achievement.'' Individuals who signed both Letters B and C, that is, signatories who both privilege debate and oppose diversity statements, are overwhelmingly inferred to be tenured white men at highly research intensive universities. Our empirical results are consistent with theories of power drawn from the social sciences.


\section*{Introduction}

The American Mathematical Society (AMS) is the largest mathematics society in the world. According to the AMS, its monthly publication, \textit{Notices of the American Mathematical Society} (\textit{AMS Notices}), is one of the most widely-read mathematics publications in the world. In the December 2019 edition of the \textit{AMS Notices}, a Vice President of the AMS published an invited essay critical of the use of diversity statements in faculty hiring within higher education \cite{Tho2019}.

In a typical faculty search process, candidates submit the following materials: a cover letter; a curriculum vitae; a research statement, which describes research experience and future plans; and a teaching statement, which discusses teaching experience and philosophy. Many faculty searches now additionally require candidates to submit a diversity statement (or similarly named document). The diversity statement is meant to empower search committees and institutions to identify candidates who have skills, experiences, and/or plans that would support inclusion, diversity, and equity on campus. For an example of an institution's rationale for and guidelines to writing such a statement, see materials from the University of California, Davis \cite{UCDavisDS}.

Ref.~\cite{Tho2019} argues that the required use of diversity statements in hiring is akin to McCarthyism, the campaign begun by Senator Joseph McCarthy against alleged communists in the United States during the early 1950s. The author writes:
\begin{quote}
In 1950 the Regents of the University of California required all UC faculty to sign a statement asserting that ``I am not a member of, nor do I support any party or organization that believes in, advocates, or teaches the overthrow of the United States Government, by force or by any illegal or unconstitutional means, that I am not a member of the Communist Party.'' Eventually thirty-one faculty members were fired over their refusal to sign... Faculty at universities across the country are facing an echo of the loyalty oath, a mandatory ``Diversity Statement'' for job applicants. The professed purpose is to identify candidates who have the skills and experience to advance institutional diversity and equity goals. In reality it's a political test, and it's a political test with teeth.
\end{quote}
The essay elicited many responses including \cite{UCDavis2019}, written by the leadership of the university where the author of \cite{Tho2019} is on faculty. In our present work, however, we focus on high-profile responses generated within the mathematical sciences community itself: one online blog post and three public letters, each with hundreds of signatories, that appeared in the \textit{AMS Notices}.

First, the Institute for the Quantitative Study of Inclusion, Diversity, and Equity published an online post \cite{QSIDE1}. This post critiques~\cite{Tho2019} and the AMS's choice to publish it. Additionally, \cite{QSIDE1} expresses concern for those in the author's professional sphere potentially affected by the publication of the essay. Finally, \cite{QSIDE1} recommends a number of actions that community members could take in response to~\cite{Tho2019}.

Second, a group of concerned mathematical scientists drafted a public letter, hereafter referred to as Letter~A \cite{letterA}. Like \cite{QSIDE1}, this letter expresses dissatisfaction with the publication of~\cite{Tho2019}. However, Letter~A highlights the need for social justice within the mathematical sciences community, assures readers that there are many individuals who see diversity work as integral, and reaffirms the importance of diversity statements in hiring.

Third, a different group of concerned mathematical scientists drafted a public letter in response to \cite{QSIDE1}, hereafter referred to as Letter~B~\cite{letterB}. This second public letter expresses displeasure with \cite{QSIDE1}, frames it as an attack on the author of \cite{Tho2019}, and provides an explicit stance neither on diversity statements nor on diversity in general. Rather, Letter~B affirms as its highest priority the need for discussion and debate within the community.

Fourth and last, another letter was drafted, hereafter referred to as Letter~C~\cite{letterC}. Unlike Letter~B, Letter~C does briefly mention the importance of ``reducing various difficulties still faced by underrepresented groups.'' However, also unlike Letter~B, it specifically addresses diversity statements, calling their mandatory use ``a mistake'' and stating that using them to eliminate candidates during the early stages of a search ``diminishes mathematical achievement.'' \nameref{S1_Appendix} contains the full text of all three public letters.

Letter~A was made available for signing beginning on November 22, 2019. Unfortunately, we do not have access to the dates that Letters~B~and~C were first made available for signing. As a result of this missing timeline information, we are not able to make meaningful statements comparing the total number of signatories to each letter. All three letters were published online on December 13, 2019, as part of the January 2020 edition of the \textit{AMS Notices} \cite{letterA,letterB,letterC}. The purpose of this manuscript is to report on crowdsourced inferences of certain demographic characteristics of the signatories.

We have declared potential conflicts of interest during the submission of this journal article, but we repeat them here in the body of our manuscript because it is an important ethical practice to do so. Co-authors CMT and JAH are co-founders of QSIDE, and CMT wrote the blog post to which Letter~B is a response. Additionally, co-authors CMT, CDE, KI, BK, DL, and JL are drafters of and signatories to Letter~A, and hence appear in the data set that we construct and study.

Our primary results are as follows. Letter~A, which highlights diversity and social justice, was signed by relatively more individuals inferred to be women and/or members of underrepresented ethnic groups. Moreover, this group is diverse with respect to the levels of professional security and types of academic institutions represented. In stark contrast, Letter~B, which privileges discussion and debate, was signed by a strong majority of individuals inferred to be white men in professionally secure positions at highly research intensive universities. Letter~C speaks out specifically against diversity statements. Individuals who signed both Letters~B~and~C, that is, signatories who privilege debate \emph{and} oppose diversity statements, are overwhelmingly inferred to be tenured white men at highly research intensive universities. These results are consistent with theories of power from the social sciences, which predict that people from demographic groups in different positions within social power structures (as determined by overrepresentation within larger populations, monetary compensation, or professional security) have divergent perspectives on these structures.

\section*{Methods}

While Letters~A~and~B were being signed electronically by the public during November and December of 2019, the lists of signatories were publicly accessible via links provided by the letters' drafters. We retrieved these lists from their respective Google spreadsheets \cite{Adata,Bdata} the morning of December 10, 2019. The AMS informed letter drafters that lists of signatories must be submitted by the evening of December 10. Therefore, the final published lists may contain a small number of additional records not included in our present study because they were submitted in the intervening hours. Additionally, for Letters~A~and~B, we removed duplicate signatures, and we removed fake signatures such as ``J. Epstein, Department of Children's Rights'' (a reference to public figure Jeffrey Epstein) and ``Donald J. Trump, University of Washington, D.C.'' It is possible that the AMS removed additional names during the publication process; our data set \emph{will} include any such names, but based on manual inspection, we believe the number to be statistically negligible. We retrieved the signatories to Letter~C once they were published \cite{letterC}. For all letters, the AMS may update the lists of signatories between the drafting of this manuscript and its publication, and our study will not reflect those changes.

To sign one of the letters, each signatory submitted their name and affiliation. We refer to this textual data as the raw data. To discover information about each signatory, we used the crowdsourcing platform Amazon Mechanical Turk (MTurk). MTurk is a crowdsourcing labor platform. Requesters sign up for the platform to post tasks that they would like completed in return for a small fee. Workers sign up for the platform and earn money by completing these tasks. For further background on MTurk, see, for example \cite{BerHubLen2012}. For our study, we served as requesters, posting each raw data record to MTurk. We asked workers to locate information about that signatory, such as personal and academic demographics. We provide details about the information requested later in this section, and we provide the full text of our Amazon Mechanical Turk questionnaire in \nameref{S1_Appendix}. A number of peer-reviewed published studies adopt this crowdsourcing methodology to infer demographics; see, for example, work on diversity in art and diversity on mathematical sciences journal editorial boards \cite{TopSen2016,TopKliTur2019}. 

We adopted several strategies to ensure data quality. First, we only hired workers who had previously completed at least 1,000 crowdworking tasks and who had a prior requester approval rating of 99\% or higher. Second, we had each signatory initially researched by four independent workers. For each item on our survey instrument, we made a final determination of the answer only if there was a consensus of at least three among the four independent responses. When we could not find a consensus, we deployed the record one additional time to MTurk and sought a consensus of three out of five. If such a consensus still did not exist, we assigned the value NA, meaning that we were not able to make an inference. For most questions, MTurk workers had the option to choose ``cannot confidently determine'' if they felt unable to make a clear inference, and we coded these responses as NA.

\subsection*{Demographic Information}

Some information we collected was for procedural use only, and is not part of our final data set. This excluded information consists of the signatory's web page, email, and year of receipt of Ph.D. (if these were available). For several reasons, our final data set does include the signatories' names. First, having appeared in the \textit{AMS Notices}, these names are public. Additionally, all information used to infer demographics is public on the internet and was accessed by MTurk workers via Wikipedia, college/university web pages, and other pages easily located from a Google search. Second, we include the signatories' names in order to follow standards of reproducible research. To omit them would make this research less verifiable. Third and last, the crowdsourcing approach to inferring demographics of public figures based on public information has been previously established in the literature \cite{TopSen2016,TopKliTur2019}. We now describe in more detail the procedures we applied to each raw data record for each question on our survey instrument that is contained in our final data set.

\subsubsection*{Affiliation} We asked workers to extract the signatory's affiliation, which was typically a company or an academic institution. In cases where multiple affiliations were listed, we accepted the first academic affiliation.

\subsubsection*{Gender} Gender is a complex construct incorporating gender identity, gender expression, social roles, and more \cite{But1990}. In reality, an individual's gender is determined by that individual. However, we judged that surveying the signatories would be unlikely to produce sufficient data, and so similar to \cite{TopSen2016,TopKliTur2019}, we have used crowdsourced inferences, recognizing the limitations of this approach, and keeping in mind that self-identifications are always preferable. We proceed with our gender analysis because both actual gender (as self-identification) and inferred gender (as perceived by others) are salient to issues of representation in the academy. We asked workers to infer the signatory's gender based on name, pictures, and/or textual information such as pronouns. The options for gender are woman, man, and nonbinary/gender nonconforming.

\subsubsection*{Ethnicity} We adopted a similar approach as for gender, again recognizing serious limitations. First, as with gender, ethnicity is most accurately stated by the individual. Second, our survey instrument asks for ``primary'' ethnicity, which disallows the identification of individuals who might be placed in multiple categories. While this approach may provide an incomplete ethnic categorization from an individual's perspective, it may reflect how an individual's ethnicity is perceived by others. We asked workers to infer ethnicity based on any information they found. The options for ethnicity are Asian, Black, Latinx, Middle Eastern, Native American/Alaska Native, Hawaiian Native/Pacific Islander, White, and Other. 

\subsubsection*{Professional sphere} We asked workers to research whether or not the signatory is situated in higher education.

\subsubsection*{Institutional Classification} If a signatory was determined to be in higher education, we asked workers to determine their type of institution, which could be a U.S. Research Intensive: Very High Research Activity (R1) institution, a U.S. Research Intensive: High Research Activity (R2) institution, another type of U.S. institution, or an institution outside of the U.S. Institutions belonging to the first two categories are determined by the Carnegie Classification of Institutions of Higher Education \cite{instlist}. We provided a list of R1 and R2 institutions to help workers make this determination. In cases where there was a consensus that the institution is outside of the U.S., we asked workers to enter the country as free text.

\subsubsection*{Professional Role} If a signatory was determined to be in higher education, we asked workers to determine their role, which could be retired/Emeritus, Professor, Associate Professor, Assistant Professor, non-tenure-track (\emph{e.g.}, lecturer, instructor, postdoctoral fellow, or visiting faculty), staff, graduate student, undergraduate student, or other. Our MTurk questionnaire comments briefly that workers should be aware of international equivalents of some of the aforementioned roles, and we check these as part of our validation process (described below).

\subsubsection*{Academic Field} If a signatory was determined to be in higher education, we asked workers to determine the person's field based on their department, degree program, title, or related information. The workers could choose Mathematics/Applied Mathematics, Statistics/Data Science, Mathematics Education, Computer Science/Computer Engineering, or other. For some individuals, multiple field associations are possible.

\subsection*{Cleaning and Validation}

After inferring demographic characteristics, part of our research team used internet searches and their own professional knowledge to fill in some missing data and to correct obvious errors, including in titles of signatories situated outside of the U.S. Additionally, after posting our data set publicly, and during the writing of this manuscript, we received emails from the community pointing out some errors in the data, namely eight for ethnicity and three for field. We made these corrections. We manually de-duplicated the data set, identifying individuals who signed more than one letter. Finally, we normalized institutional affiliations so that, for instance, ``UC Davis'' and ``University of California Davis'' would be aggregated together.

Then, a separate part of our research team manually checked the accuracy of a random 5\% subsample of the records having a professional affiliation in higher education. For the 72 records checked, we agreed with all crowdsourced inferences for gender, institutional classification, and academic field. For inferred ethnicity, we detected one possible error for an individual who was categorized as White by MTurk workers, but whom we ourselves would have classified as ``cannot confidently determine.'' This accounts for an observed 1.4\% error rate for inferred ethnicity in our data set. However, a 95\% confidence interval for the true error rate reaches approximately 4\%, meaning (loosely) that amongst the 1,367 higher education signatories in our data set, we might reasonably expect up to 55 errors in inferred ethnicity. We detected four errors in professional role. Two individuals, labeled respectively as Assistant Professor and Associate Professor, should have both been categorized as non-tenure-track. Two individuals labeled as Professor should have been labeled as retired/Emeritus. Later, we will aggregate professional roles into two coarser categories of professional security. The four aforementioned errors only result in one error under that coarser categorization: an Associate Teaching Professor classified as more professionally secure should have been classified as less professionally secure because the position is not tenurable. This is the same error rate as for inferred ethnicity. We found no other errors. While we cannot know the true error rates in our data, our validation procedure suggests that they are small, especially for our coarser categorizations. After completing our validation procedure, we edited our data set to correct the aforementioned errors we found.

\section*{Results}

Letter~A has 621 signatories, Letter~B has 680 signatories, and Letter~C has 210 signatories. The overlap between signatories to Letter~A and Letter~B is six individuals. There is no overlap between the signatories to Letters~A~and~C. There is substantial overlap between the signatories to Letters~B~and~C, namely 74 individuals. The six individuals who signed Letters~A~and~B do not constitute a group large enough to analyze meaningfully, so we leave them in that data set accompanies this manuscript, but we exclude them from our analyses (meaning they are not represented in the tables and charts prsented here). We divide the remaining signatories into four disjoint groups: the 615 signatories to Letter~A only, the 600 signatories to Letter~B only, the 136 signatories to Letter~C only, and the 74 signatories to Letters~B~and~C, for a total of 1,425 signatories. Hereafter, we refer to the groups as Pool~A~Only, Pool~B~Only, Pool~C~Only, and Pool~B~and~C.

We divide our presentation of results into two subsections. First, we examine personal demographics, namely inferred gender and inferred ethnicity. Second, for signatories situated in higher education, which constitute over 95\% of the original data, we examine institutional classification, professional role, and academic field. For each of the aforementioned variables except for inferred gender, we also develop a simplified version of the variable by aggregating certain levels of the original one. We describe these aggregations later.

For each demographic variable, we present a two-way contingency table showing proportions of letter signatory pool versus that variable. See Table~\ref{table1} for personal demographics and Table~\ref{table2} for academic demographics. Within a given contingency table, for each combination of pool and demographic, the left-hand number gives the column percentage and the right-hand number gives the row percentage. For example, for inferred gender in Table~\ref{table1}, looking down the first column, Pool~A~Only has 0.2\% signatories for whom we did not infer gender, 0.2\% of individuals inferred as nonbinary, 50.7\% inferred as women, and 48.9\% inferred as men. In the same section of Table~\ref{table1}, looking across the row tabulating women, we see that of all individuals inferred to be women in our data set, 69.2\% of them are in Pool~A~Only, 21.1\% are in Pool~B~Only, 8.2\% are in Pool~C~Only, and 1.6\% are in Pool~B~and~C.

We perform statistical tests on the contingency tables for gender and for the simplified variables (described below) with NA values excluded. For these tables, we first perform a $\chi^2$ test of the null hypothesis that the demographic characteristic and letter signatory pool are independent. At the top of the contingency table, we report the value of the $\chi^2$ test statistic, the degrees of freedom $df$, the $p$-value, and Cram\'{e}r's effect size $V$. All $\chi^2$ tests result in rejection of the null hypothesis. Second, we test for significance of individual cells in the contingency table via a $z$-test for sample percentage. Colors indicate whether, at the $\alpha = 0.05$ significance level, a particular cell represents more (blue) or fewer (orange) individuals than we would expect if the variables were independent. We use Holm-adjusted $p$-values to correct for multiple comparisons. For background on the statistical approaches we use, see reference texts such as \cite{Agr2018}.

It is important to remember that the population in our statistical tests is the 1,425 letter signatories. As we highlight later, these signatories do not represent a random sample of the mathematical sciences community (or any other group). In the presentation of our results, we will use the words ``overrepresented'' and ``underrepresented.'' When we use these words in the following two subsections, we mean over/under-represented as compared to what we would expect if demographics and signatory pool were independent.

\begin{table}[!t]
\small
\centering
\begin{adjustwidth}{-2.25in}{0in}
\caption{
{\bf Two-way contingency tables for letter signatory pool and personal demographics, namely inferred gender, inferred ethnicity, and inferred underrepresented minority (URM) status. Gender and URM both fail a statistical test for the independence of demographics and pool.} Within a given contingency table, for each combination of pool and demographic, the left-hand number gives the column percentage and the right-hand number gives the row percentage. Colors indicate whether a particular cell represents more (blue) or fewer (orange) individuals than we would expect if the variables were independent. For details of statistical tests, see table notes below.}
\begin{tabular}{lcccccccc}
\thickhline
\textbf{Pool}                                     & \multicolumn{2}{c}{\textbf{A Only}}                                                                     & \multicolumn{2}{c}{\textbf{B Only}}                                                                     & \multicolumn{2}{c}{\textbf{C Only}}                             & \multicolumn{2}{c}{\textbf{B and C}}                                                                  \\
                                                  & \multicolumn{2}{c}{\textbf{n = 615}}                                                                    & \multicolumn{2}{c}{\textbf{n = 600}}                                                                    & \multicolumn{2}{c}{\textbf{n = 136}}                            & \multicolumn{2}{c}{\textbf{n = 74}}                                                                   \\
                                                  & Col \%                                             & Row \%                                             & Col \%                                             & Row \%                                             & Col \%                         & Row \%                         & Col \%                                            & Row \%                                            \\ \thickhline
\textbf{Inferred Gender}                          &                                                    &                                                    &                                                    &                                                    &                                &                                &                                                   &                                                   \\
$\chi^2 = 189.8$, $df=3$, $p<0.00001$, $V=0.366$  & \multicolumn{1}{l}{}                               & \multicolumn{1}{l}{}                               & \multicolumn{1}{l}{}                               & \multicolumn{1}{l}{}                               & \multicolumn{1}{l}{}           & \multicolumn{1}{l}{}           & \multicolumn{1}{l}{}                              & \multicolumn{1}{l}{}                              \\
NA (n = 7)                                        & 0.2\%                                              & 14.3\%                                             & 1.0\%                                              & 85.7\%                                             & 0.0\%                          & 0.0\%                          & 0.0\%                                             & 0.0\%                                             \\
Nonbinary (n = 1)                                 & 0.2\%                                              & 100.0\%                                            & 0.0\%                                              & 0.0\%                                              & 0.0\%                          & 0.0\%                          & 0.0\%                                             & 0.0\%                                             \\
Woman (n = 451)                                   & \multicolumn{1}{l}{\cellcolor[HTML]{9698ED}50.7\%} & \multicolumn{1}{l}{\cellcolor[HTML]{9698ED}69.2\%} & \multicolumn{1}{l}{\cellcolor[HTML]{FE996B}15.8\%} & \multicolumn{1}{l}{\cellcolor[HTML]{FE996B}21.1\%} & \multicolumn{1}{l}{27.2\%}     & \multicolumn{1}{l}{8.2\%}      & \multicolumn{1}{l}{\cellcolor[HTML]{FE996B}9.5\%} & \multicolumn{1}{l}{\cellcolor[HTML]{FE996B}1.6\%} \\
Man (n = 966)                                     & \cellcolor[HTML]{FE996B}48.9\%                     & \cellcolor[HTML]{FE996B}31.2\%                     & \cellcolor[HTML]{9698ED}83.2\%                     & \cellcolor[HTML]{9698ED}51.7\%                     & 72.8\%                         & 10.2\%                         & \cellcolor[HTML]{9698ED}90.5\%                    & \cellcolor[HTML]{9698ED}6.9\%                     \\ \thickhline
\textbf{Inferred Ethnicity}                       &                                                    &                                                    &                                                    &                                                    &                                &                                &                                                   &                                                   \\
statistical tests not performed                   & \multicolumn{1}{l}{}                               & \multicolumn{1}{l}{}                               & \multicolumn{1}{l}{}                               & \multicolumn{1}{l}{}                               & \multicolumn{1}{l}{}           & \multicolumn{1}{l}{}           & \multicolumn{1}{l}{}                              & \multicolumn{1}{l}{}                              \\
NA (n = 80)                                       & 6.3\%                                              & 48.8\%                                             & 6.3\%                                              & 47.5\%                                             & 2.2\%                          & 3.8\%                          & 0.0\%                                             & 0.0\%                                             \\
Asian (n = 108)                                   & 7.0\%                                              & 39.8\%                                             & 7.8\%                                              & 43.5\%                                             & 13.2\%                         & 16.7\%                         & 0.0\%                                             & 0.0\%                                             \\
Black (n = 27)                                    & 3.3\%                                              & 74.1\%                                             & 1.0\%                                              & 22.2\%                                             & 0.7\%                          & 3.7\%                          & 0.0\%                                             & 0.0\%                                             \\
Latinx (n = 68)                                   & 8.5\%                                              & 76.5\%                                             & 2.2\%                                              & 19.1\%                                             & 0.7\%                          & 1.5\%                          & 2.7\%                                             & 2.9\%                                             \\
ME (n = 42)                                       & 2.1\%                                              & 31.0\%                                             & 4.2\%                                              & 59.5\%                                             & 2.2\%                          & 7.1\%                          & 1.4\%                                             & 2.4\%                                             \\
NAAN (n = 1)                                      & 0.2\%                                              & 100.0\%                                            & 0.0\%                                              & 0.0\%                                              & 0.0\%                          & 0.0\%                          & 0.0\%                                             & 0.0\%                                             \\
HNPI (n = 2)                                      & 0.3\%                                              & 100.0\%                                            & 0.0\%                                              & 0.0\%                                              & 0.0\%                          & 0.0\%                          & 0.0\%                                             & 0.0\%                                             \\
White (n = 1097)                                  & 72.4\%                                             & 40.6\%                                             & 78.5\%                                             & 42.9\%                                             & 80.9\%                         & 10.0\%                         & 95.9\%                                            & 6.5\%                                             \\ \hline
\textbf{Inferred Underrepresented Status}         &                                                    &                                                    &                                                    &                                                    &                                &                                &                                                   &                                                   \\
$\chi^2$ = 49.6, $df=3$, $p < 0.00001$, $V=0.192$ & \multicolumn{1}{l}{}                               & \multicolumn{1}{l}{}                               & \multicolumn{1}{l}{}                               & \multicolumn{1}{l}{}                               & \multicolumn{1}{l}{}           & \multicolumn{1}{l}{}           & \multicolumn{1}{l}{}                              & \multicolumn{1}{l}{}                              \\
NA (n = 80)                                       & 6.3\%                                              & 48.8\%                                             & 6.3\%                                              & 47.5\%                                             & 2.2\%                          & 3.8\%                          & 0.0\%                                             & 0.0\%                                             \\
Not URM (n = 1247)                                & \cellcolor[HTML]{FE996B}81.5\%                     & \cellcolor[HTML]{FE996B}40.2\%                     & \cellcolor[HTML]{9698ED}90.5\%                     & \cellcolor[HTML]{9698ED}43.5\%                     & \cellcolor[HTML]{9698ED}96.7\% & \cellcolor[HTML]{9698ED}10.5\% & 97.3\%                                            & 5.8\%                                             \\
URM (n = 98)                                      & \cellcolor[HTML]{9698ED}12.2\%                     & \cellcolor[HTML]{9698ED}76.5\%                     & \cellcolor[HTML]{FE996B}3.2\%                      & \cellcolor[HTML]{FE996B}19.4\%                     & \cellcolor[HTML]{FE996B}1.5\%  & \cellcolor[HTML]{FE996B}2.0\%  & 2.7\%                                             & 2.0\%                                             \\ \thickhline
\end{tabular}
\begin{flushleft} Results of a $\chi^2$ test for independence appear at the top of the contingency tables for inferred gender and inferred underrepresented minority status. Here, $\chi^2$ is the value of the test statistic, $df$ is degrees of freedom, $p$ is the probability value, and $V$ is Cram\'{e}r's effect size. We perform the $\chi^2$ tests on modified contingency tables (not shown) with the small number of NA values removed. To test for significance of individual cells (as coded by the colors) we use a $z$-test for sample percentage at the $\alpha=0.05$ significance level with Holm-adjusted $p$-values to account for multiple comparisons. Inferred underrepresented status is a collapsed form of the inferred ethnicity variable. Above, we use the following abbreviations: ME (Middle Eastern), NAAN (Native American/Alaska Native), HNPI (Hawaiian Native/Pacific Islander) and URM (ethnic groups underrepresented in mathematics, namely Black, Latinx, NAAN, and HNPI). Non-URM ethnic groups are White, Middle Eastern, and Asian.
\end{flushleft}
\label{table1}
\end{adjustwidth}
\end{table}

\begin{table}[!ht]
\small
\centering
\begin{adjustwidth}{-2.25in}{0in}
\caption{
{\bf Two-way contingency tables for letter signatory pool and academic demographics, namely institutional classification, research intensiveness, professional role, professional security, academic field, and simplified academic field. Research intensiveness, professional security, and simplified field each fails a statistical test for the independence of demographics and pool.} Within a given contingency table, for each combination of pool and demographic, the left-hand number gives the column percentage and the right-hand number gives the row percentage. Colors indicate whether a particular cell represents more (blue) or fewer (orange) individuals than we would expect if the variables were independent. For details of statistical tests, see table notes below.}
\begin{tabular}{lcccccccc}
\thickhline
\textbf{Pool}                                      & \multicolumn{2}{c}{\textbf{A~Only~HE}}                          & \multicolumn{2}{c}{\textbf{B~Only~HE}}                                                                        & \multicolumn{2}{c}{\textbf{C~Only~HE}}                                                                      & \multicolumn{2}{c}{\textbf{B~and~C~HE}}                                                                     \\
                                                   & \multicolumn{2}{c}{\textbf{n = 586}}                            & \multicolumn{2}{c}{\textbf{n = 563}}                                                                          & \multicolumn{2}{c}{\textbf{n = 135}}                                                                        & \multicolumn{2}{c}{\textbf{n = 73}}                                                                         \\
                                                   & Col \%                         & Row \%                         & Col \%                                                & Row \%                                                & Col \%                                               & Row \%                                               & Col \%                                               & Row \%                                               \\ \thickhline
\textbf{Institutional Classification}              &                                &                                &                                                       &                                                       &                                                      &                                                      &                                                      &                                                      \\
statistical tests not performed                    & \multicolumn{1}{l}{}           & \multicolumn{1}{l}{}           & \multicolumn{1}{l}{}                                  & \multicolumn{1}{l}{}                                  & \multicolumn{1}{l}{}                                 & \multicolumn{1}{l}{}                                 & \multicolumn{1}{l}{}                                 & \multicolumn{1}{l}{}                                 \\
US R1 Very High Research Activity (n = 672)        & 37.2\%                         & 32.4\%                         & 52.4\%                                                & 43.9\%                                                & 71.1\%                                               & 14.3\%                                               & 86.3\%                                               & 9.4\%                                                \\
US R2 High Research Activity (n = 103)             & 11.1\%                         & 63.1\%                         & 6.2\%                                                 & 34.0\%                                                & 2.2\%                                                & 2.9\%                                                & 0.0\%                                                & 0.0\%                                                \\
US Other (n = 377)                                 & 46.2\%                         & 71.9\%                         & 18.1\%                                                & 27.1\%                                                & 2.2\%                                                & 0.8\%                                                & 1.4\%                                                & 0.3\%                                                \\
Non US (n = 205)                                   & 5.5\%                          & 15.6\%                         & 23.3\%                                                & 63.9\%                                                & 24.4\%                                               & 16.1\%                                               & 12.3\%                                               & 4.4\%                                                \\ \hline
\textbf{Research Intensiveness}                    &                                &                                &                                                       &                                                       &                                                      &                                                      &                                                      &                                                      \\
$\chi^2 = 247.1$, $df=3$, $p < 0.00001$, $V=0.427$ & \multicolumn{1}{l}{}           & \multicolumn{1}{l}{}           & \multicolumn{1}{l}{}                                  & \multicolumn{1}{l}{}                                  & \multicolumn{1}{l}{}                                 & \multicolumn{1}{l}{}                                 & \multicolumn{1}{l}{}                                 & \multicolumn{1}{l}{}                                 \\
Less Research Intensive (n = 480)                  & \cellcolor[HTML]{9698ED}57.3\% & \cellcolor[HTML]{9698ED}70.0\% & \cellcolor[HTML]{FE996B}24.3\%                        & \cellcolor[HTML]{FE996B}28.5\%                        & \cellcolor[HTML]{FE996B}4.4\%                        & \cellcolor[HTML]{FE996B}1.2\%                        & \cellcolor[HTML]{FE996B}1.4\%                        & \cellcolor[HTML]{FE996B}0.2\%                        \\
More Research Intensive (n = 877)                  & \cellcolor[HTML]{FE996B}42.7\% & \cellcolor[HTML]{FE996B}28.5\% & \cellcolor[HTML]{9698ED}75.7\%                        & \cellcolor[HTML]{9698ED}48.6\%                        & \cellcolor[HTML]{9698ED}95.6\%                       & \cellcolor[HTML]{9698ED}14.7\%                       & \cellcolor[HTML]{9698ED}98.6\%                       & \cellcolor[HTML]{9698ED}8.2\%                        \\ \thickhline
\textbf{Professional Role}                         &                                &                                &                                                       &                                                       &                                                      &                                                      &                                                      &                                                      \\
statistical tests not performed                    & \multicolumn{1}{l}{}           & \multicolumn{1}{l}{}           & \multicolumn{1}{l}{}                                  & \multicolumn{1}{l}{}                                  & \multicolumn{1}{l}{}                                 & \multicolumn{1}{l}{}                                 & \multicolumn{1}{l}{}                                 & \multicolumn{1}{l}{}                                 \\
NA (n = 34)                                        & 3.1\%                          & 52.9\%                         & 2.7\%                                                 & 44.1\%                                                & 0.0\%                                                & 0.0\%                                                & 1.4\%                                                & 2.9\%                                                \\
Undergraduate Student (n = 7)                      & 0.5\%                          & 42.9\%                         & 0.7\%                                                 & 57.1\%                                                & 0.0\%                                                & 0.0\%                                                & 0.0\%                                                & 0.0\%                                                \\
Graduate Student (n = 119)                         & 14.8\%                         & 73.1\%                         & 5.7\%                                                 & 26.9\%                                                & 0.0\%                                                & 0.0\%                                                & 0.0\%                                                & 0.0\%                                                \\
Staff (n = 11)                                     & 1.2\%                          & 63.6\%                         & 0.7\%                                                 & 36.4\%                                                & 0.0\%                                                & 0.0\%                                                & 0.0\%                                                & 0.0\%                                                \\
Non-Tenure-Track Position (n = 89)                 & 10.1\%                         & 66.3\%                         & 5.0\%                                                 & 31.5\%                                                & 0.7\%                                                & 1.1\%                                                & 1.4\%                                                & 1.1\%                                                \\
Assistant Professor (n = 194)                      & 24.9\%                         & 75.3\%                         & 8.5\%                                                 & 24.7\%                                                & 0.0\%                                                & 0.0\%                                                & 0.0\%                                                & 0.0\%                                                \\
Associate Professor (n = 207)                      & 20.5\%                         & 58.0\%                         & 13.9\%                                                & 37.7\%                                                & 5.9\%                                                & 3.9\%                                                & 1.4\%                                                & 0.5\%                                                \\
Full Professor (n = 629)                           & 24.1\%                         & 22.4\%                         & 54.2\%                                                & 48.5\%                                                & 85.9\%                                               & 18.4\%                                               & 91.8\%                                               & 10.7\%                                               \\
Retired/Emeritus Professor (n = 67)                & 0.9\%                          & 7.5\%                          & 8.7\%                                                 & 73.1\%                                                & 7.4\%                                                & 14.9\%                                               & 4.1\%                                                & 4.5\%                                                \\ \hline
\textbf{Professional Security}                     &                                &                                &                                                       &                                                       &                                                      &                                                      &                                                      &                                                      \\
$\chi^2$ = 238.2, $df=3$, $p < 0.00001$, $V=0.424$ & \multicolumn{1}{l}{}           & \multicolumn{1}{l}{}           & \multicolumn{1}{l}{}                                  & \multicolumn{1}{l}{}                                  & \multicolumn{1}{l}{}                                 & \multicolumn{1}{l}{}                                 & \multicolumn{1}{l}{}                                 & \multicolumn{1}{l}{}                                 \\
NA (n = 35)                                        & 3.2\%                          & 54.3\%                         & 2.7\%                                                 & 42.9\%                                                & 0.0\%                                                & 0.0\%                                                & 1.4\%                                                & 2.9\%                                                \\
Less Secure (n = 419)                              & \cellcolor[HTML]{9698ED}51.4\% & \cellcolor[HTML]{9698ED}71.8\% & \cellcolor[HTML]{FE996B}{\color[HTML]{000000} 20.6\%} & \cellcolor[HTML]{FE996B}{\color[HTML]{000000} 27.7\%} & \cellcolor[HTML]{FE996B}{\color[HTML]{000000} 0.7\%} & \cellcolor[HTML]{FE996B}{\color[HTML]{000000} 0.2\%} & \cellcolor[HTML]{FE996B}{\color[HTML]{000000} 1.4\%} & \cellcolor[HTML]{FE996B}{\color[HTML]{000000} 0.2\%} \\
More Secure (n = 903)                              & \cellcolor[HTML]{FE996B}45.4\% & \cellcolor[HTML]{FE996B}29.5\% & \cellcolor[HTML]{9698ED}76.7\%                        & \cellcolor[HTML]{9698ED}47.8\%                        & \cellcolor[HTML]{9698ED}99.3\%                       & \cellcolor[HTML]{9698ED}14.8\%                       & \cellcolor[HTML]{9698ED}97.3\%                       & \cellcolor[HTML]{9698ED}7.9\%                        \\ \thickhline
\textbf{Field}                                     &                                &                                &                                                       &                                                       &                                                      &                                                      &                                                      &                                                      \\
statistical tests not performed                    & \multicolumn{1}{l}{}           & \multicolumn{1}{l}{}           & \multicolumn{1}{l}{}                                  & \multicolumn{1}{l}{}                                  & \multicolumn{1}{l}{}                                 & \multicolumn{1}{l}{}                                 & \multicolumn{1}{l}{}                                 & \multicolumn{1}{l}{}                                 \\
NA (n = 13)                                        & 0.9\%                          & 38.5\%                         & 1.4\%                                                 & 61.5\%                                                & 0.0\%                                                & 0.0\%                                                & 0.0\%                                                & 0.0\%                                                \\
Math/Applied Math (n = 1206)                       & 87.5\%                         & 42.5\%                         & 87.0\%                                                & 40.6\%                                                & 96.3\%                                               & 10.8\%                                               & 100.0\%                                              & 6.1\%                                                \\
Statistics/Data Science (n = 12)                   & 1.5\%                          & 75.0\%                         & 0.5\%                                                 & 25.0\%                                                & 0.0\%                                                & 0.0\%                                                & 0.0\%                                                & 0.0\%                                                \\
Mathematics Education (n = 36)                     & 6.0\%                          & 97.2\%                         & 0.2\%                                                 & 2.8\%                                                 & 0.0\%                                                & 0.0\%                                                & 0.0\%                                                & 0.0\%                                                \\
Computer Science (n = 16)                          & 1.2\%                          & 43.8\%                         & 1.2\%                                                 & 43.8\%                                                & 1.5\%                                                & 12.5\%                                               & 0.0\%                                                & 0.0\%                                                \\
Other (n = 74)                                     & 2.9\%                          & 23.0\%                         & 9.6\%                                                 & 73.0\%                                                & 2.2\%                                                & 4.1\%                                                & 0.0\%                                                & 0.0\%                                                \\ \hline
\textbf{Simplified Field}                          &                                &                                &                                                       &                                                       &                                                      &                                                      &                                                      &                                                      \\
$\chi^2$ = 67.9, $df=6$, $p < 0.00001$, $V=0.159$  & \multicolumn{1}{l}{}           & \multicolumn{1}{l}{}           & \multicolumn{1}{l}{}                                  & \multicolumn{1}{l}{}                                  & \multicolumn{1}{l}{}                                 & \multicolumn{1}{l}{}                                 & \multicolumn{1}{l}{}                                 & \multicolumn{1}{l}{}                                 \\
NA (n = 13)                                        & 0.9\%                          & 38.5\%                         & 1.4\%                                                 & 61.5\%                                                & 0.0\%                                                & 0.0\%                                                & 0.0\%                                                & 0.0\%                                                \\
MASD (n = 1218)                                    & 89.1\%                         & 42.9\%                         & 87.6\%                                                & 40.5\%                                                & 96.3\%                                               & 10.7\%                                               & \cellcolor[HTML]{9698ED}100.0\%                      & \cellcolor[HTML]{9698ED}6.0\%                        \\
Mathematics Education (n = 36)                     & \cellcolor[HTML]{9698ED}6.0\%  & \cellcolor[HTML]{9698ED}97.2\% & \cellcolor[HTML]{FE996B}0.2\%                         & \cellcolor[HTML]{FE996B}2.8\%                         & 0.0\%                                                & 0.0\%                                                & 0.0\%                                                & 0.0\%                                                \\
Computer Science/Other (n = 90)                    & \cellcolor[HTML]{FE996B}4.1\%  & \cellcolor[HTML]{FE996B}26.7\% & \cellcolor[HTML]{9698ED}10.8\%                        & \cellcolor[HTML]{9698ED}67.8\%                        & 3.7\%                                                & 5.6\%                                                & 0.0\%                                                & 0.0\%                                                \\ \thickhline
\end{tabular}
\begin{flushleft} Results of a $\chi^2$ test for independence appear at the top of the contingency tables for research intensiveness, professional security, and simplified field. Here, $\chi^2$ is the value of the test statistic, $df$ is degrees of freedom, $p$ is the probability value, and $V$ is Cram\'{e}r's effect size. When NA values are present, we perform the $\chi^2$ tests on modified contingency tables (not shown) with the small number of NA values removed. To test for significance of individual cells (as coded by the colors) we use a $z$-test for sample percentage at the $\alpha=0.05$ significance level with Holm-adjusted $p$-values to account for multiple comparisons. Research intensiveness, professional security, and simplified field are collapsed versions of institutional classification, professional role, and field. Above, we use the abbreviation MASD (Mathematics/Applied Mathematics/Statistics/Data Science). 
\end{flushleft}
\label{table2}
\end{adjustwidth}
\end{table}

\subsection*{Personal Demographics}

We reiterate that all demographic data are \emph{inferred} data. This is especially critical for gender and ethnicity. Below, we will use language such as ``individuals inferred to be women'' and ``individuals identified as belonging to underrepresented ethnic groups.'' Ideally, gender and ethnicity would be self-identified.

Fig~\ref{fig1}(A) displays the inferred gender composition of each letter signatory pool; the percentages shown are column percentages from Table~\ref{table1}. The proportion of individuals inferred to be nonbinary or gender nonconforming is low. We do not conjecture on the true proportion of these individuals within our data, nor within the mathematical sciences community at large. However, we do suspect that the low percentage we observe is attributable at least in part to the limitations of our method of gender inference \cite{Nic2019}. As for the remaining (binary) gender categories, Pool~A~Only achieves approximate gender parity, comprising 50.7\% inferred women. In contrast, Pools~B Only, C~Only, and B~and~C have low percentages -- 15.8\%, 27.2\%, and 9.5\% respectively -- and hence are dominated by individuals inferred to be men. Compared to what we would expect under the assumption of independence of gender and pool, signatories identified as women are overrepresented in Pool~A~Only and underrepresented in Pool~B~Only and Pool~B~and~C.

\begin{figure}[t!]
\caption{{\bf Personal demographic characteristics of letter signatory pools.} See Table~\ref{table1} for details. (A)~Inferred gender. Pool~A~Only achieves (approximate) gender parity, whereas Pools~B~Only, C~Only, and B~and~C are dominated by people inferred to be men. (B)~Inferred underrepresented minority (URM) status. At 12.2\%, Pool~A~Only has a higher than expected percentage of people inferred to belong to ethnic groups traditionally underrepresented in the mathematical sciences. Pools~B~Only, C~Only, and B~and~C each have approximately 3\% or less.}
\label{fig1}
\end{figure}

Fig~\ref{fig1}(B) examines inferred ethnicity. Table~\ref{table1} shows the full results for our original ethnicity variable. However, because there are many ethnicity categories, we perform an additional analysis by grouping together minoritized ethnicities that are considered to be traditionally underrepresented (URM) in mathematics: Black, Latinx, Native American/Alaska Native (NAAN), and Hawaiian Native/Pacific Islander (HNPI). Under this classification, the percentage inferred to be URM is 12.2\% for Pool~A~Only, 3.2\% for Pool~B~Only, 1.5\% for Pool~C~Only, and 2.7\% for Pool~B~and~C, as shown in Fig~\ref{fig1}(B). The statistical tests of Table~\ref{table1} indicate that individuals identified as URM are overrepresented in Pool~A~Only and underrepresented in Pools~B~Only and C~Only.

\subsection*{Academic Demographics}

We now restrict attention to the signatories situated within higher education (HE). We refer to the pools as Pool~A~Only~HE, Pool~B~Only~HE, Pool~C~Only~HE, and Pool~B~and~C~HE. These pools have, respectively, 586, 563, 135, and 73 members, for a total of 1,357 individuals.

Fig~\ref{fig2}(A) addresses institutional classification. Table~\ref{table2} shows the full results for our original institutional classification variable. We perform an additional analysis by creating a collapsed version of this variable as we now describe. The non-U.S. countries appearing most frequently in our data (and their counts) are Israel (47), Canada (37), the United Kingdom (21), Germany (14), and France (13). The specific non-U.S. institutions appearing most frequently in our data (and their counts) are Hebrew University of Jerusalem (20), Technion (15), and University of Toronto (13). We judge that if located in the U.S., these institutions would be classified as Very High Research Activity (R1). Therefore, we aggregate US R1 and international institutions into the category of more research intensive (more RI), and we aggregate US R2 and US other institutions into the category of less research intensive (less RI). As shown in Fig~\ref{fig2}(A), Pool~A~Only~HE is composed 57.3\% of less RI, whereas Pools~B~Only~HE, C~Only~HE, and B~and~C~HE have 24.3\%, 4.4\%, and 1.4\% respectively. The latter two pools are especially dominated by individuals from more RI institutions, whereas Pool~A~Only HE is fairly balanced. Indeed, our statistical tests indicate that less RI institutions are underrepresented in the last three pools.

\begin{figure}[ht!]
\caption{{\bf Academic demographic characteristics of letter signatory pools.} See Table~\ref{table2} for details. (A)~Research intensiveness. Pool~A~Only~HE shows a balance of signatories. Pools~B~Only~HE, C~Only~HE, and especially B~and~C~HE are dominated by people from highly research intensive institutions. (B)~Professional security. At 51.4\%, Pool~A~Only~HE overrepresents individuals with less professional security. Pools~B~Only~HE, C~Only~HE, and B~and~C~HE underrepresent these individuals, and the latter two contain nearly zero less professionally secure individuals. (C)~Academic field. Pool~A~Only~HE contains the highest percentage of signatories from mathematics education, who are overrepresented at 6.0\%.}
\label{fig2}
\end{figure}

Fig~\ref{fig2}(B) addresses professional security. Table~\ref{table2} shows the full results for our original professional role variable. We perform an additional analysis by creating a collapsed version of this variable. We aggregate the roles that are not eligible for academic tenure or are eligible but have not yet received tenure (as judged based on title), and refer to these as less professionally secure roles. These roles are undergraduate student, graduate student, staff, non-tenure-track faculty, and Assistant Professor (and their international equivalents). Additionally, we aggregate the roles that do indicate currently or previously held academic tenure, namely Associate Professor, Full Professor, and retired/Emeritus (and their international equivalents). We refer to these as more professionally secure roles. As shown in Fig~\ref{fig2}(B), in Pool~A~Only~HE, more professionally secure individuals are in the minority, at 45.4\%.  Pools~B~Only~HE, C~Only~HE, and B~and~C~HE are dominated by more professionally secure individuals, at 76.7\%, 99.3\%, and 97.3\% respectively. More professionally secure individuals are statistically overrepresented in these pools. Notably, if we examine Pools~C~Only~HE and B~and~C~HE, we see that almost every signatory to Letter~C is tenured (or was at some point).

Finally, Fig~\ref{fig2}(C) addresses academic field. Table~\ref{table2} shows the full results for our original academic field variable. We perform an additional analysis by creating a collapsed version of this variable. We aggregate the mathematics/applied mathematics and statistics/data science categories into one category, and we aggregate computer science/computer engineering and other disciplines into another category, leaving mathematics education as the third category. As shown in Fig~\ref{fig2}(C), Pool~A~Only~HE contains the highest percentage of signatories from mathematics education, who are statistically overrepresented at 6.0\%. Pool~B~Only~HE contains a substantial group of signatories from computing and other fields, who are statistically overrepresented at 10.8\%. Pool~B~and~C~HE is composed exclusively of mathematicians and applied mathematicians.

\subsection*{Comparison with the Profession}

To contextualize our results, we compare gender and ethnicity percentages in the mathematical sciences at large to those of letter signatories in higher education whose field is identified as mathematics, applied mathematics, statistics, or data science. When we refer to demographic percentages below, we always mean ``within the mathematical sciences.'' The set of comparisons we are able to make is limited by the availability of data on the profession.

\subsubsection*{Gender}

Women account for approximately 40\% of undergraduate mathematical sciences degrees \cite{AMSsurvey1}. There are only seven undergraduate students in our data set, and we could not infer gender for one of them. The remaining six are inferred as three women in Pool~A~Only~HE and three men in Pool~B~Only~HE.

At the doctoral level, women account for approximately 30\% of mathematical sciences degrees awarded \cite{AMSsurvey2}. Our data set contains 118 graduate students for whom we inferred gender, and none signed Letter~C. In Pool~A~Only~HE, individuals inferred to be women are overrepresented compared to the field at large, comprising 46.1\%. In Pool~B~Only~HE, they are underrepresented, comprising only 17.2\%.

Women comprise approximately 20\% of tenure-stream faculty at doctoral degree granting departments of mathematical sciences in the United States \cite{AMSsurvey3}. To make a comparison with our study, we restrict attention to Assistant Professors, Associate Professors, and Professors at U.S. R1 and R2 institutions, of which there are 499 in our data set. Individuals inferred to be women are overrepresented in Pool~A~Only~HE, comprising 46.9\%. They are underrepresented in Pool~B~Only~HE, comprising 10.2\%. In Pool~C~Only~HE, they comprise 23.0\%, just over the at-large value. In Pool~B~and~C~HE, they are underrepresented, comprising 11.9\%.

\subsubsection*{Ethnicity}

There is not good availability of data on ethnicity in the profession. For instance, the percentage of active faculty having URM status is unknown \cite{QSIDEfoc}. We do know that individuals having URM status comprise 8\% of the pool of doctoral degrees granted in the mathematical sciences \cite{AMSsurvey2}. In our data set, there are 91 graduate students for whom we inferred ethnicity. None of them signed Letter~C. In Pool~A~Only~HE, individuals inferred as URM are overrepresented, comprising 18.8\%, more than double the at-large value. In Pool~B~Only~HE, they comprise 9.1\%, approximately on par with the at-large value.

\section*{Discussion and Conclusions}

We have conducted a crowdsourced study of the demographics of signatories to three public letters in the mathematical sciences community. Our key results are as follows. As for personal demographics, we infer Pool~A~Only to be substantially more diverse than the other pools. We have calculated for each pool the (joint) percentage of signatories who were inferred to be men and who were not inferred to be members of an underrepresented ethnic group. For Pool~A~Only, this percentage is 40.8\%, whereas for Pools~B~Only, C~Only, and B~and~C, the percentage is 75.8\%, 70.6\%, and 87.8\% respectively. As for academic demographics, Pool~A~Only represents a broader range of institution types and levels of professional security. We have calculated for each pool the (joint) percentage of signatories who are situated at highly research intensive institutions in more professionally secure positions. For Pool~A~Only~HE, this percentage is merely 11.8\%, whereas for Pools~B~Only~HE, C~Only~HE, and B~and~C~HE, the percentage is 55.6\%, 94.8\%, and 95.9\% respectively. Finally, restricting attention to our HE pools, we have calculated the percentages of individuals who (jointly) are inferred to be men, not members of underrepresented ethnic groups, located at highly research intensive institutions, and in more professionally secure positions. For Pool~A~Only~HE, this percentage is a scant 5.6\%. For the remaining pools it is 46.5\%, 67.4\%, and 83.6\%.

In summary, Letter~A highlights diversity and social justice and was signed by relatively more people inferred to be women and members of minoritized ethnic groups. These signatories represent a broader range of institution types and levels of professional security. Letter~B does not comment on diversity, but rather, argues for discussion and debate. It was signed predominantly by individuals inferred to be men who have ethnicities not underrepresented in the mathematical sciences, and who are in professionally secure positions at highly research intensive institutions. Letter~C speaks out specifically against the use of diversity statements. Individuals who signed both Letters~B~and~C, that is, signatories who privilege discussion and debate and who oppose diversity statements, are overwhelmingly inferred to be tenured white men at highly research intensive universities.

We now relate our results to theories of power drawn from scholarship in the social sciences. Our findings are consistent with the idea of \emph{positionality}. Positionality describes the ways in which individuals' identities and experiences are consequences of their positions within social structures. These positions shape an individual's perceptions \cite{Alc1988}. More specifically, positionality theory predicts that individuals' positions within power structures tilt their perceptions of phenomena in patterned ways \cite{Har1987}. Differences in perception are particularly pronounced when individuals' identities confer unequal levels of power. For example, among students who were shown the same film about race relations, white students were more likely to respond by describing the film as an exaggeration, or ``over the top,'' while students of color were more likely to respond with reflections on how much the film mirrored their own experiences \cite{KinOrb2008}. In this example, students' own racial identifications (their positions within social structures) produced very different perceptions about the same film (the phenomenon).

In our present study, positionality theory suggests that people with relatively more power in the mathematical sciences would have very different perspectives from those with relatively less power. This accords with the aggregate patterns we have documented. People in more powerful positions within the mathematical sciences, namely men, white people, people with tenure, and people at highly research intensive universities, tend in aggregate to endorse perspectives on diversity statements that are different from the perspectives endorsed in aggregate by people with relatively less power, namely women, members of underrepresented ethnic groups, people without tenure, and people at institutions that are not highly research intensive.

Some might read tenure and affiliation with highly research intensive universities as proxies for the quality of a mathematician. We are not persuaded that this reading undermines our interpretation above. Regardless of how validly or invalidly tenure status and institution type capture the quality of a mathematician, these professional attributes are face-valid indicators of power. Both tenure and employment at highly research intensive institutions confer substantial monetary benefits. The wage premium associated with highly research intensive universities is especially marked in the natural and mathematical sciences \cite{MelStr2007}. This premium persists even if one controls for individuals' research productivity, years of experience, and demographics. All other things being equal, being tenured and affiliated with highly research intensive institutions means more access to social power than not being tenured, or being employed at other types of institutions. Moreover, tenure confers increased academic freedom \cite{Che1998} and professional security. Money, academic freedom, and professional security all help individuals exercise their will, or equivalently, give individuals more power \cite{Web1978}.  

Additionally, our results are consistent with prior studies that use discourse analysis. When individuals denounce actions that they perceive as harmful to a group, members of other groups often discursively frame this denunciation as an attempt to silence the critiqued actors \cite{Dij1992}. Prior research has documented this dynamic in contexts including standup comedy \cite{Per2013}, college campus controversies \cite{MooBel2017},
the immigration policy debate \cite{Chi2010}, and more. The framing of denunciation as an attempt at silencing is consistent with the statements in Letter B, as endorsed by its signatories; see \nameref{S1_Appendix}. Indeed, this letter begins with the sentence ``We write with grave concerns about recent attempts to intimidate a voice within our mathematical community.''

In our study, we have reported descriptive statistics and conducted selected statistical tests. Future work could undertake further statistical modeling and inference. Finally, while results of our study are consistent with theories of power, our study is not explanatory. Though Pool~A~Only has markedly different demographics from the other pools, we do not know why this is case. For instance, demographic differences could arise from the way that news of the letters was disseminated through professional networks, reflecting or even amplifying the (potential) demographic biases of those networks. Alternatively, demographic differences between the pools could reflect contrasting personal and professional values. Professional security could also have played a role in a signing decision of Letters B and C. Of course, these demographic differences could be caused by a combination of the aforementioned factors and others. Investigation using tools and frameworks from the social sciences and humanities might complement our research to provide explanations and further insights.

\section*{Supporting information}


\label{S1_Appendix}

\subsection*{S1 Appendix A: Letter A}

We are a group of concerned mathematicians writing in response to AMS Vice President Abigail Thompson's editorial, invited by the AMS for publication in the December 2019 edition of the Notices. In this editorial, Dr. Thompson states her personal opinion against the mandated use of faculty diversity statements in hiring decisions and compares such requirements to McCarthyist loyalty oaths.

We are all members of many mathematical societies, including the American Mathematical Society. Some of us serve on committees in these societies or are chairs of committees in these societies. Some of us are or have been chairs of departments, some of us are or have been chairs of search committees, and some of us have written or reviewed diversity statements as part of search processes. We have all thought deeply about the role of diversity statements and related tools, such as student success statements.

We are compelled to write because the AMS leadership's actions have harmed the mathematics community, particularly mathematicians from marginalized backgrounds. We are writing because we support diversity statements as one tool to encourage a more inclusive and equitable mathematics profession. We are writing because we wish to correct the misleading impressions readers might have of such statements from Thompson's editorial: Thompson's opinion does not represent the opinions of many other members of the mathematics community. We are writing because not everyone is in a position to raise their voice. We are writing because it matters how our community and its leaders talk about diversity, especially in our profession's most prominent publication. We are writing because we are disappointed by the editorial decision to publish the piece which contradicts the AMS's commitment to diversity affirmed in its own diversity statement. Clearly, this is something that people needed to talk about, but the AMS has chosen to spark this conversation by giving its imprimatur to a piece that undermines productive discussion and causes real danger and burden to the marginalized members of our community.

Diversity statements are widely used in academic hiring as one component to assess candidates' qualifications for the job. Each statement one requires as part of a hiring process – research, teaching, mentoring, service, or diversity – helps paint a picture of how a candidate will contribute to the work of an institution. Increased use of diversity statements reflects a growing recognition in higher education that faculty contribute in positive ways to the campus community by acknowledging, appreciating, and collaborating with groups of students, staff, and fellow faculty who are diverse along varied axes. In acknowledgement that this is part of the work of a faculty member and of the hiring process, we recommend that graduate programs explicitly prepare their graduates to contribute to this work and to write and talk about it meaningfully, and we commend the programs already undertaking this work. There are plenty of legitimate questions about how to use diversity statements effectively and how (more broadly) to create diverse and supportive faculties. In order to reduce bias in the evaluation of candidates, hiring committees evaluate statements according to criteria that indicate evidence of these important contributions, grounded in the missions of higher education in general and their institution in particular. Asking for and evaluating diversity statements are not quick solutions to the complex challenge of justice and inclusion in higher education, but they can help hiring committees to evaluate candidates' skills in doing this portion of our professional work.

Diversity statements help assess a candidate's ability to effectively teach a diverse group of students. If our goal as mathematicians and educators is truly to reach as many students as possible, thinking about diversity and inclusion is necessary. Good teaching is necessarily inclusive. If we willfully ignore an important area of pedagogy that demonstrably helps more students succeed in math, then we will continue to reproduce systems of inequity, and we will do a great disservice to our students. We will therefore not be effective teachers.

Suggesting that actively attempting to include more students in mathematics is equivalent to the Red Scare is ignorant (about both history and the present) and dangerous. Claims of ``reverse racism,'' which equate critiquing privilege with oppressing the privileged, have a long and unsavory history in and beyond higher education. Without understanding the history in which these discussions are rooted, it is possible to profess support for the ideal of equality while acting in ways that lead to exclusion and inequity.

While Dr. Thompson attempts to spin this issue with partisan wording, diversity statements are a small yet necessary step towards creating a more equitable and inclusive community. Higher education in the US is shifting, student populations we serve are changing, and our understanding of how to better serve all students is advancing. We need a rehumanization of mathematics that can affirm students' cultural funds of knowledge while examining and combating its own roles in supporting power structures. We need leadership at all levels, from professional societies to presidents, boards, deans, and chairs, to recognize this reality, advocate for students and faculty from a variety of backgrounds, and move us forward.

Dr. Thompson's preface that the letter is her ``personal opinion'' does not alleviate our concerns, nor does the fact that she seems to be referring primarily to the use of these documents at the UC system. The fact remains that the Notices made an editorial decision to give Thompson's essay a national (indeed, international) platform, and in a prominent position within the publication. Notices is a publication of the AMS, and Dr. Thompson is identified as an AMS officer in her byline. According to Notices editor Erica Flapan, Dr. Thompson's position in the AMS leadership led the AMS to solicit her letter. These contextual details send a message to the profession about how diversity is viewed by those with power and responsibility in the AMS and a major university department. The AMS and Notices bear responsibility for amplifying views that fly in the face of research-based practices and that falsely equate evidence-based approaches to teaching and professional practice with the blacklisting of people based on political ideology, all in direct contradiction of the AMS's stated commitment to diversity.

AMS's own diversity statement claims, ``The American Mathematical Society is committed to promoting and facilitating equity, diversity and inclusion throughout the mathematical sciences… We reaffirm the pledge in the AMS Mission Statement to `advance the status of the profession of mathematics, encouraging and facilitating full participation of all individuals,' and urge all members to conduct their professional activities with this goal in mind.'' While merely publishing Dr. Thompson's letter demonstrates the AMS's lack of commitment to this statement, the fact that it was written by and credited to an officer of the AMS raises even more serious questions about the statement's sincerity.

We strongly disagree with the sentiments and arguments in Dr. Thompson's editorial, and we hope that the AMS will reconsider the way that it uses its power and position in the mathematics communities in these kinds of discussions. However, we primarily write this letter to our fellow mathematicians and students of all kinds who might have wondered if inclusion work is valued in our community. We want students and faculty, especially those with multiple identities that are minoritized in mathematics, to know that many mathematicians see this inclusion work as integral to our community and identity. 

\subsection*{S1 Appendix B: Letter B}

We write with grave concerns about recent attempts to intimidate a voice within our mathematical community. Abigail Thompson published an opinion piece in the December issue of the Notices of the American Mathematical Society. She explained her support for efforts within our community to further diversity, and then described her concerns with the rigid rubrics used to evaluate diversity statements in the hiring processes of the University of California system.

The reaction to the article has been swift and vehement. An article posted at the site QSIDE urges faculty to direct their students not to attend and not to apply for jobs at the University of California, Davis, where Prof. Thompson is chair of the math department. It recommends contacting the university to question whether Prof. Thompson is fit to be chair. And it recommends refusing to do work for the Notices of the American Mathematical Society for allowing this piece to be published.

Regardless of where anyone stands on the issue of whether diversity statements are a fair or effective means to further diversity aims, we should agree that this attempt to silence opinions is damaging to the profession. This is a direct attempt to destroy Prof. Thompson's career and to punish her department. It is an attempt to intimidate the AMS into publishing only articles that hew to a very specific point of view. If we allow ourselves to be intimidated into avoiding discussion of how best to achieve diversity, we undermine our attempts to achieve it.

We the undersigned urge the American Mathematical Society to stand by the principle that important issues should be openly discussed in a respectful manner, and to make a clear statement that bullying and intimidation have no place in our community.

\subsection*{S1 Appendix C: Letter C}

In an essay in the December 2019 issue of the Notices, Abigail Thompson describes the mandatory ``Diversity Statement'' (mDS) that mathematics job applicants to UC Schools must submit together with their regular applications. At some campuses, the mDS is evaluated, in various categories, according to a detailed list of criteria (called a ``rubric'').

If an applicant to Berkeley, for example, merely says that she advocates ``mentoring,  treating  all  students  the  same  regardless of background,'' she merits a score of 1--2 out of a possible top score of 5 in the ``track record for advancing diversity''  category*.  Hiring  committees  (at  UC Davis and Berkeley, in particular) are encouraged by the  Administration  to  use  the  rubrics,  establish  a  cutoff,  and  eliminate  candidates  who  score  below  the  cutoff  as  a  first  step  in  the  hiring  process  for  all  hires.  In  this  way  Diversity Statements diminish the value of mathematical achievement as the key element in securing a position at a UC Mathematics Department.

Mandatory  Diversity  Statements  undermine  Faculty  Governance. Should the use of scored diversity statements become required as the first step in the hiring process, this opens  the  way  for  Administrators,  who  have  no  professional knowledge of mathematics, to exert primary control over  the  hiring  of  mathematicians.  And  indeed,  testing  the  waters,  small  scale  pilot  programs  have  already  been  implemented at various UC schools requiring the first cut on hiring to be based on such scored diversity statements.

We applaud Abigail Thompson for her courageous leadership in bringing this issue to the attention of the broader Mathematics Community. As she says in her essay:

``Mathematics has made progress over the past decades towards becoming a more welcoming, inclusive discipline.  We should continue to do all we can to reduce barriers to participation  in  this  most  beautiful  of  fields....  There  are  reasonable means to further this goal: encouraging students from all backgrounds to enter the mathematics pipeline, trying to ensure that talented mathematicians don't leave the profession, creating family-friendly policies, and sup-porting junior faculty at the beginning of their careers, for example.''

We  agree  wholeheartedly  with  these  sentiments.  It  is  important  to  strive  to  hire  faculty  that  will  make  the  atmosphere more welcoming to all. It is also important to recognize  and  help  reduce  various  difficulties  still  faced  by underrepresented groups. But as Abigail says, there are mistakes to avoid: mDS's are one of them.

Finally, we commend the Editorial Board at the Notices for  opening  up  the  discussion  on  this  very  important  matter.

*If you insert the following into a google search, the first thing which comes up is a pdf with the rubric for Berkeley:
\medskip
\\
\small\noindent{\texttt{rubric\_to\_assess\_candidate\_contributions\_to\_diversity\_equity\_and\_inclusion-1.pdf}}

\subsection*{S1 Appendix D: Amazon Mechanical Turk Questionnaire}

We are interested in gathering demographic information about people who signed a high profile public letter to an academic/professional society. Immediately below, you will see the way a person signed the letter.

\medskip
\noindent {[public signature text information]}

\medskip
\noindent We ask you to find certain demographic information about the person above. To find this information, you might use:
\begin{itemize}
\item    Wikipedia
\item    Google searches
\item    The Mathematics Genealogy Project
\item    any other resources you find helpful
\end{itemize}

\medskip
\noindent 1. The information above needs to be separated into a person's name and professional affiliation. Please enter the person's name as they have written it above. That is, separate the person's name from all other information, and enter it here. If there are typos or extra spaces, please clean these up.

\medskip
\noindent {[free text entry]}

\medskip
\noindent 2. The information above needs to be separated into a person's name and professional affiliation. Please enter the person's professional affiliation. That is, separate the affiliation from all other information and enter it here. This should be the name of an educational institution or a company ONLY. If the person has given other information (for instance, title, academic department, professional specialization, honors/awards, random geographic info that is not part of the company's or academic institution's name) please ignore it in answering this question. If the person has listed multiple affiliations, enter only the first one. If the person has not listed an affiliation (just a name) then see if you can locate the person's affiliation (by Googling) and enter it here. If you cannot determine their affiliation with confidence, please leave this question blank.

\medskip
\noindent {[free text entry]}

\medskip
\noindent 3. Please try to locate web information associated with the person. This might be their own personalized web site and/or the web site of a department, college/university, or company. Please enter at most one website, even if you use more in answering other questions below. If you cannot locate any web information, please leave this question blank.

\medskip
\noindent {[free text entry]}

\medskip
\noindent 4. Please try to locate an email address associated with the person. A professional email address (college/university, company) is preferred. If you cannot locate an email address, please leave this question blank.

\medskip
\noindent {[free text entry]}

\medskip
\noindent 5. What do you think the person's gender is? In answering this question, you might use their name, and/or any pictures you find of them, and/or any textual information you find referring to them with gendered pronouns.
\begin{itemize}
\item    Woman
\item    Man
\item    Binary man/woman choice does not apply (nonbinary, gender nonconconforming)
\item     Cannot confidently determine (insufficient time/information)
\end{itemize}

\medskip
\noindent 6. An ethnic group is a category of people who identify with each other based on similarities such as common ancestral, linguistic, social, cultural or national experiences. What is the person's primary ethnic background?
\begin{itemize}
\item     American Indian or Alaska Native. A person having origins in any of the original peoples of North and South America (including Central America), and who maintains tribal affiliation or community attachment.
 \item    Black or African American. A person having origins in any of the black racial groups of Africa.
\item     Asian. A person having origins in any of the original peoples of the Far East, Southeast Asia, or the Indian subcontinent including, for example, Cambodia, China, India, Japan, Korea, Malaysia, Pakistan, the Philippine Islands, Thailand, and Vietnam.
 \item    Latinx. A person of Cuban, Mexican, Puerto Rican, South or Central American origin, or similar.
 \item    Native Hawaiian or Other Pacific Islander. A person having origins in any of the original peoples of Hawaii, Guam, Samoa, or other Pacific Islands.
\item     Middle East or North Africa. A person having origins in locations such as Algeria, Morocco, Egypt, Lebanon, or Syria.
\item     White. A person having origins in Europe.
\item     Other (not listed above).
 \item    Cannot confidently determine (insufficient time/information).
 \end{itemize}
 
\medskip
\noindent 7. Based on the person's affiliation as given above (and, as necessary, based on other information you find) does the person have an affiliation with higher education (student or employee at a college, university, or other post-secondary educational institution)? Note: please answer yes for a person who is retired or ``Emeritus'' from teaching in higher education.
\begin{itemize}
\item     Yes, the person is affiliated with higher education.
 \item    No, the person is not affiliated with higher education.
 \item    Cannot confidently determine (insufficient time/information).
\end{itemize}

\medskip
\noindent \textbf{If the answer to question 7 above is NO or CANNOT DETERMINE, please scroll down to the bottom of this form and submit it without answering any further questions. If the answer to question 7 is YES, please continue.}

\medskip
\noindent 8. We would like to know a little bit more about the person's college/university. In particular we'd like to know if it is in the United States, and if so, its degree of research intensiveness. As for research intensiveness, please see this Wikipedia page:

\medskip
\noindent \url{https://en.wikipedia.org/wiki/List\_of\_research\_universities\_in\_the\_United\_States}

\medskip
\noindent A college/university may be in the first list on that page as ``R1: Doctoral Universities – Very high research activity,'' or in the second list on that page as ``R2: Doctoral Universities – High research activity,'' or on neither list. Choose the appropriate option below.
\begin{itemize}
\item     The institution is in the United States and is on the R1 list.
\item     The institution is in the United States and is on the R2 list.
\item     The institution is in the United States but is on neither the R1 list nor the R2 list.
\item     The institution is not in the United States.
 \item    Cannot confidently determine (insufficient time/information).
 \end{itemize}
 
 \medskip
\noindent 9. If you selected NOT IN THE UNITED STATES in question 8 above, please type the full name of the country where the institution is located. Otherwise, please leave this question blank.

\medskip
\noindent {[free text entry]}

\medskip
\noindent 10. We would like to know a little bit more about the person's role/title at their college/university. Choose the appropriate option below.
\begin{itemize}
\item     The person is a retired professor (they may or may not have the title Emeritus).
\item     Their title is Professor (or an international equivalent).
 \item    Their title is Associate Professor (or an international equivalent).
 \item    Their title is Assistant Professor (or an international equivalent).
 \item    The person is an instructor who is not eligible for academic tenure. For instance, their title is Instructor or Lecturer or Postdoctoral Fellow, or it contains one of the following phrases: ``NTT'', ``non-tenure track,'' ``visiting'', or similar.
\item     The person is a graduate student (earning master's degree or PhD).
\item     The person is an undergraduate student (earning bachelor's degree or similar).
 \item    The person is staff at the institution (they are not a student nor are they any of the faculty types listed above).
 \item    Other (none of the choices above).
 \item    Cannot confidently determine (insufficient time/information).
 \end{itemize}

\medskip
\noindent 11. If you chose STAFF, OTHER, or CANNOT DETERMINE for question 10, please skip this question. Otherwise, we would like to know a little bit more about the person's academic field. This might be based on their title, academic department, degree program, or other information. Choose the appropriate option below.
\begin{itemize}
\item     Mathematics or Applied Mathematics.
 \item    Statistics or Data Science.
 \item    Computer Science, Computer Engineering, or similar.
 \item    Mathematics Education (this means that the person's research and/or studies pertain to how best to teach mathematics).
 \item    Other (none of the choices above).
 \item    Cannot confidently determine (insufficient time/information).
 \end{itemize}

\medskip
\noindent 12. If in question 10 you chose GRADUATE STUDENT, UNDERGRADUATE STUDENT, STAFF, OTHER, or CANNOT DETERMINE, please skip this question. Otherwise, for question 8, you indicated one of the first five choices (some type of professor or other instructor). Determine the year that person earned a PhD degree. You might find this on their web page, their curriculum vitae, or in the Mathematics Genealogy Project. Please enter the Ph.D. year as a four-digit number. If you cannot find this information, or if you believe the person does not have a Ph.D., please leave this question blank.

\medskip
\noindent{[free text entry]}

\section*{Acknowledgments}
We are grateful to Alicia Prieto Langarica for contributions during the early stages of this work, and to Carina Curto and Joshua Paik for helpful comments regarding our methodology and data.

\nolinenumbers

\bibliography{AMSlettersbib}




\end{document}